\documentclass[12pt,twoside]{article}
\usepackage[english]{babel}
\usepackage[latin1]{inputenc}
\usepackage{amsmath}
\usepackage{amssymb,amsfonts}
\usepackage{graphicx}
\usepackage{times,amssymb,amscd}

\newcommand{\bC}{\mathbf{C}}
\newcommand{\bE}{\mathbf{E}}
\newcommand{\bG}{\mathbf{G}}
\newcommand{\bH}{\mathbf{H}}

\newcommand{\bL}{\mathbf{L}}

\newcommand{\bR}{\mathbf{R}}
\newcommand{\bS}{\mathbf{S}}

\newcommand{\bs}{\mathbf{s}}
\newcommand{\ba}{\mathbf{a}}

\newcommand{\bT}{\mathbf{T}}

\newcommand{\bb}{\mathbf{b}}

\newcommand{\cP}{\mathcal{P}}
\newcommand{\cS}{\mathcal{S}}
\newcommand{\cT}{\mathcal{T}}

\newcommand{\cH}{\mathcal{H}}

\newcommand{\EUC}{\bE^3}

\newcommand{\HYP}{\bH^3}
\newcommand{\SXR}{\bS^2\!\times\!\bR}
\newcommand{\HXR}{\bH^2\!\times\!\bR}
\newcommand{\SLR}{\widetilde{\bS\bL_2\bR}}
\newcommand{\NIL}{\mathbf{Nil}}
\newcommand{\SOL}{\mathbf{Sol}}

\newtheorem{Definition}{Definition}[section]
\newtheorem{Remark}{Remark}[section]

\begin{document}
\pagestyle{myheadings}
\markboth{\centerline{Emil Moln\'ar and Jen\H o Szirmai}}
{Volumes and geodesic ball packings to the regular prism tilings...}
\title
{Volumes and geodesic ball packings to the regular prism tilings in $\SLR$ space \footnote{Mathematics Subject Classification 2010: 52C17, 52C22, 52B15, 53A35, 51M20. \newline
Key words and phrases: Thurston geometries, $\SLR$ geometry, density of ball packing under space group, regular
prism tiling, volume in the $\SLR$ space.}}

\vspace{2mm}
\author{\vspace{3mm} \normalsize{Dedicated to Prof. Lajos Tamássy on His $90^{th}$ Birthday}\\
Emil Moln\'ar and Jen\H o Szirmai \\
\normalsize Budapest University of Technology and \\
\normalsize Economics Institute of Mathematics, \\
\normalsize Department of Geometry \\
\normalsize Budapest, P. O. Box: 91, H-1521 \\
\normalsize emolnar@math.bme.hu,~szirmai@math.bme.hu
\date{\normalsize{\today}}}

%%%%%%%%%%%%%%%%%%%%%%%%%%%%%%%%%%%%%%%%%%%%
%%AMS Classification 2000
%% The 1-st classification is obligatory, the 2-nd classification is optional
%% \subjclass{primary}{secondary}      f.e. \subjclass{35R35, 49N50}{}
%%{52C17, 52C22}{}
%%%%%%%%%%%%%%%%%%%%%%%%%%%%%%%%%%%%%%%%%%%%

\maketitle
\begin{abstract}

After having investigated the regular prisms and prism tilings in the $\SLR$ space in the previous work \cite{Sz13-1} of the second author, we consider
the problem of geodesic ball packings related to those tilings and their symmetry groups $\mathbf{pq2_1}$. $\SLR$ is one of
the eight Thurston geometries that can be derived from the 3-dimensional Lie group of all $2\times 2$ real matrices with determinant one.

In this paper we consider geodesic spheres and balls in $\SLR$ (even in $\mathbf{SL_{\mathrm{2}}R})$, if their radii $\rho\in [0, \frac{\pi}{2})$,
and determine their volumes. Moreover, we consider the prisms of the above space and compute their volumes,
define the notion of the geodesic ball packing and its density. We develop a procedure to determine
the densities of the densest geodesic ball packings for the tilings, or in this paper more precisely, for their generating groups $\mathbf{pq2_1}$
(for integer rotational parameters $p,q$; $3\le p,~\frac{2p}{p-2} <q$). 
We look for those parameters $p$ and $q$ above, where the packing density large enough as possible. Now our record is $0.567362$ for $(p, q) = (8, 10)$.
These computations seem to be important, since we do not know optimal ball packing, namely in the hyperbolic space $\HYP$.
We know only the density upper bound 0.85326, realized by horoball packing of $\HYP$ to its ideal regular simplex tiling.  
Surprisingly, for the so-called translation ball packings under the same groups $\mathbf{pq2_1}$ in \cite{MSzV13} 
we have got larger density $0.841700$ for $(p, q) = (5, 10000 \rightarrow \infty)$ close to the above upper bound.

We use for the computation and visualization of the $\SLR$ space its projective model introduced by the first author in \cite{M97}.

\end{abstract}
%%%%%%%%%%%%%%%%%%%%%%%%%%%%%%%%%%%%%%%%%%%%

%%%%%%%%%%%%%%%%%%%%%%%%%%%%%%%%%%%%%%%%%%%
\newtheorem{theorem}{Theorem}[section]
\newtheorem{corollary}[theorem]{Corollary}
\newtheorem{conjecture}[theorem]{Conjecture}
\newtheorem{lemma}[theorem]{Lemma}
\newtheorem{exmple}[theorem]{Example}
\newenvironment{definition}{\begin{defn}\normalfont}{\end{defn}}
\newenvironment{remark}{\begin{rmrk}\normalfont}{\end{rmrk}}
\newenvironment{example}{\begin{exmple}\normalfont}{\end{exmple}}
\newenvironment{acknowledgement}{Acknowledgement}

%%%%%%%%%%%%%%%%%%%%%%%%%%%%%%%%%%%%%%%%%%%%%%%%%%%%%%%%%%%%%%%%%%%%

%============================================================================%
%                             the main article                               %
%============================================================================%

%%%%%%%%%%%%%%%%%%%%%%%%%%%%%%%%%%%%%%%%%%%%%%%%%%%%%%%%%%%%%%%%%%%%%%%%%%%%%%
\section{On $\SLR$ geometry}

The real $ 2\times 2$ matrices $\begin{pmatrix}
         d&b \\
         c&a \\
         \end{pmatrix}$ with unit determinant $ad-bc=1$
constitute a Lie transformation group by the usual product operation, taken to act on row matrices as on point coordinates on the right as follows
\begin{equation}
\begin{gathered}
(z^0,z^1)\begin{pmatrix}
         d&b \\
         c&a \\
         \end{pmatrix}=(z^0d+z^1c, z^0 b+z^1a)=(w^0,w^1) \\ 
\mathrm{with} \ w=\frac{w^1}{w^0}=\frac{b+\frac{z^1}{z^0}a}{d+\frac{z^1}{z^0}c}=\frac{b+za}{d+zc} \tag{1.1}
\end{gathered}
\end{equation}
as right action on the complex projective line $\bC^\infty$.
This group is a $3$-dimensional manifold, because of its $3$ independent real coordinates and with its usual neighbourhood topology
\cite{S}, \cite{T}.
In order to model the above structure in the projective sphere $\cP \cS^3$ and in the projective space $\cP^3$ (see \cite{M97}),
we introduce the new projective coordinates $(x^0,x^1,x^2,x^3)$ where
\begin{equation}
a:=x^0+x^3, \ b:=x^1+x^2, \ c:=-x^1+x^2, \ d:=x^0-x^3, \notag
\end{equation}
with positive, then the non-zero multiplicative equivalence as a projective freedom in $\cP \cS^3$ and in $\cP^3$, respectively. Meanwhile we turn
to the proportionality $\mathbf{SL_{\mathrm{2}}R} < \mathbf{PSL_{\mathrm{2}}R}$, natural in this context.
Then it follows that
\begin{equation}
0>bc-ad=-x^0x^0-x^1x^1+x^2x^2+x^3x^3 \tag{1.2}
\end{equation}
describes the interior of the above one-sheeted hyperboloid solid $\cH$ in the usual Euclidean coordinate simplex, with the origin
$E_0(1;0;0;0)$ and the ideal points of the axes $E_1^\infty(0;1;0;0)$, $E_2^\infty(0;0;1;0)$, $E_3^\infty(0;0;0;1)$.
We consider the collineation group ${\bf G}_*$ that acts on the projective sphere $\cS \cP^3$  and preserves a polarity, i.e. a scalar product of signature
$(- - + +)$, this group leaves the one sheeted hyperboloid solid $\cH$ invariant.
We have to choose an appropriate subgroup $\mathbf{G}$ of $\mathbf{G}_*$ as isometry group, then the universal covering group and space
$\widetilde{\cH}$ of $\cH$ will be the hyperboloid model of $\SLR$ (see Fig.~1 and \cite{M97}).

The specific isometries $\bS(\phi)$ $(\phi \in \bR )$ constitute a one parameter group given by the matrices
\begin{equation}
\begin{gathered} \bS(\phi):~(s_i^j(\phi))=
\begin{pmatrix}
\cos{\phi}&\sin{\phi}&0&0 \\
-\sin{\phi}&\cos{\phi}&0&0 \\
0&0&\cos{\phi}&-\sin{\phi} \\
0&0&\sin{\phi}&\cos{\phi}
\end{pmatrix}.
\end{gathered} \tag{1.3}
\end{equation}
The elements of $\bS(\phi)$ are the so-called {\it fibre translations}. We obtain a unique fibre line to each $X(x^0;x^1;x^2;x^3) \in \widetilde{\cH}$
as the orbit by right action of $\bS(\phi)$ on $X$. The coordinates of points lying on the fibre line through $X$ can be expressed
as the images of $X$ by $\bS(\phi)$:
\begin{equation}
\begin{gathered}
(x^0;x^1;x^2;x^3) \stackrel{\bS(\phi)}{\longrightarrow} {(x^0 \cos{\phi}-x^1 \sin{\phi}; x^0 \sin{\phi} + x^1 \cos{\phi};} \\ {x^2 \cos{\phi} + x^3 \sin{\phi};-x^2 \sin{\phi}+
x^3 \cos{\phi})}
\end{gathered} \tag{1.4}
\end{equation}
for the Euclidean coordinates $x:=\frac{x^1}{x^0}$, $y:=\frac{x^2}{x^0}$, $z:=\frac{x^3}{x^0}$, $x^0\ne 0$ as well.
The $\pi$ periodicity for the above coordinates in the above maps can be seen from the formula (1.4).
\begin{figure}[ht]
\centering
\includegraphics[width=8cm]{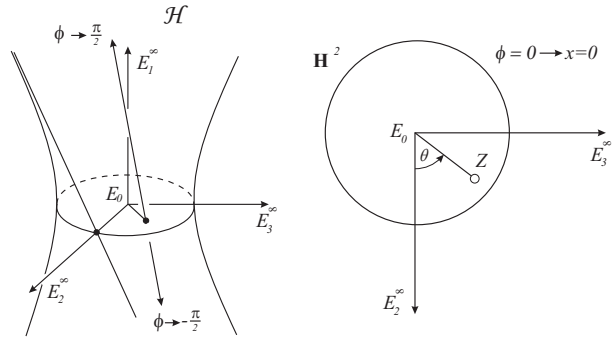}
\caption{The hyperboloid model}
\label{}
\end{figure}
In (1.3) and (1.4) we can see the $2\pi$ periodicity of $\phi$. Moreover, we see the (logical) extension to $\phi \in \bR$, as real parameter, to have
the universal covers $\widetilde{\cH}$ and $\SLR$, respectively, through the projective sphere $\cP\cS^3$. The elements of the isometry group of
$\mathbf{SL_{\mathrm{2}}R}$ (and so by the above extension the isometries of $\SLR$) can be described by the matrix $(a_i^j)$ (see \cite{M97} and \cite{MSz})
\begin{equation}
\begin{gathered} (a_i^j)=
\begin{pmatrix}
a_0^0&a_0^1&a_0^2&a_0^3 \\
\mp a_0^1&\pm a_0^0&\pm a_0^3&\mp a_0^2 \\
a_2^0&a_2^1&a_2^2&a_2^3 \\
\pm a_2^1&\mp a_2^0&\mp a_2^3&\pm a_2^2 \\
\end{pmatrix} \ \ \text{where} \\
-(a_0^0)^2-(a_0^1)^2+(a_0^2)^2+(a_0^3)^2=-1, \ \ -(a_2^0)^2-(a_2^1)^2+(a_2^2)^2+(a_2^3)^2=1, \\
-a_0^0a_2^0-a_0^1a_2^1+a_0^2a_2^2+a_0^3a_2^3=0=-a_0^0a_2^1+a_0^1a_2^0-a_0^2a_2^3+a_0^3a_2^2,
\end{gathered} \tag{1.5}
\end{equation}
and we allow positive proportionality, of course, as projective freedom.
We define the {\it translation group} $\bG_T$, as a subgroup of the isometry group of $\mathbf{SL_{\mathrm{2}}R}$,
those isometries acting transitively on the points of ${\cH}$ and by the above extension on the points of $\widetilde{\cH}$.
$\bG_T$ maps the origin $E_0(1;0;0;0)$ onto $X(x^0;x^1;x^2;x^3)$. These isometries and their inverses (
up to a positive determinant factor) can be given by 
\begin{equation}
\begin{gathered} \bT:~(t_i^j)=
\begin{pmatrix}
x^0&x^1&x^2&x^3 \\
-x^1&x^0&x^3&-x^2 \\
x^2&x^3&x^0&x^1 \\
x^3&-x^2&-x^1&x^0
\end{pmatrix},\\ \bT^{-1}:~(T_j^k)=
\begin{pmatrix}
x^0&-x^1&-x^2&-x^3 \\
x^1&x^0&-x^3&x^2 \\
-x^2&-x^3&x^0&-x^1 \\
-x^3&x^2&x^1&x^0
\end{pmatrix}.
\end{gathered} \tag{1.6}
\end{equation}
The rotation about the fibre line through the origin $E_0(1;0;0;0)$ by angle $\omega$ $(-\pi<\omega\le \pi)$ can be expressed by
\begin{equation}
\begin{gathered} \bR_{E_O}(\omega):~(r_i^j(E_0,\omega))=
\begin{pmatrix}
1&0&0&0 \\
0&1&0&0 \\
0&0&\cos{\omega}&\sin{\omega} \\
0&0&-\sin{\omega}&\cos{\omega}
\end{pmatrix},
\end{gathered} \tag{1.7}
\end{equation}
and the rotation $\bR_X(\omega)$ with matrix $:~(r_i^j(X,\omega))$ about the fibre line through $X(x^0;x^1;x^2;x^3)$ by 
angle $\omega$ can be derived by formulas (1.6) and (1.7) by conjugacy $\bR_X(\omega)=\bT^{-1} \bR_{E_O} (\omega) \bT$.
Thus the above rotation $\bR_X(\omega)$, with a specific $X$ $(\cosh{r},0,\sinh{r},0) \sim (1,0,\tanh{r},0)$ has the important matrix (see \cite{Sz13-1})
{\small{\begin{equation}
\begin{pmatrix}
\begin{gathered} 1+\sinh^2{r}- \\ -\sinh^2{r} \cos{\omega} \end{gathered} & \sinh^2{r}\sin{\omega}& \begin{gathered} \frac{1}{2}\sinh{2r}- \\
-\frac{1}{2} \sinh{2r} \cos{\omega} \end{gathered} &-\frac{1}{2}\sinh{2r}\sin{\omega} \\
-\sinh^2{r}\sin{\omega}&\begin{gathered} 1+\sinh^2{r}- \\ -\sinh^2{r} \cos{\omega} \end{gathered} &-\frac{1}{2}\sinh{2r}\sin{\omega}&\begin{gathered}
-\frac{1}{2}\sinh{2r}+ \\
+\frac{1}{2} \sinh{2r} \cos{\omega} \end{gathered} \\
\begin{gathered}
-\frac{1}{2}\sinh{2r}+ \\
+\frac{1}{2} \sinh{2r} \cos{\omega} \end{gathered} &-\frac{1}{2}\sinh{2r}\sin{\omega}&\begin{gathered} 1-\cosh^2{r}+ \\ +\cosh^2{r} \cos{\omega}
\end{gathered}
&\cosh^2{r}\sin{\omega} \\
-\frac{1}{2}\sinh{2r}\sin{\omega} &\begin{gathered}
\frac{1}{2}\sinh{2r}- \\
-\frac{1}{2} \sinh{2r} \cos{\omega} \end{gathered}&-\cosh^2{r}\sin{\omega} &\begin{gathered} 1-\cosh^2{r}+ \\ +\cosh^2{r} \cos{\omega} \end{gathered}
\end{pmatrix} \tag{1.8}
\end{equation}}}
Horizontal intersection of the hyperboloid solid $\cH$  with the plane $E_0 E_2^\infty E_3^\infty$
provides the {\it base plane} of the model $\widetilde{\cH}=\SLR$.
The fibre through $X$ intersects the hyperbolic $(\mathbf{H}^2)$ base plane $z^1=x=0$ in the foot point
\begin{equation}
\begin{gathered}
Z(z^0=x^0 x^0+x^1x^1; z^1=0; z^2=x^0x^2-x^1x^3;z^3=x^0x^3+x^1x^2).
\end{gathered} \tag{1.9}
\end{equation}
We generally introduce a so-called hyperboloid parametrization by \cite{M97} as follows
\begin{equation}
\begin{gathered}
x^0=\cosh{r} \cos{\phi}, \\
x^1=\cosh{r} \sin{\phi}, \\
x^2=\sinh{r} \cos{(\theta-\phi)}, \\
x^3=\sinh{r} \sin{(\theta-\phi)},
\end{gathered} \tag{1.10}
\end{equation}
where $(r,\theta)$ are the polar coordinates of the $\mathbf{H}^2$ base plane, and $\phi$ is the fibre coordinate. We note that
$$-x^0x^0-x^1x^1+x^2x^2+x^3x^3=-\cosh^2{r}+\sinh^2{r}=-1<0.$$
The inhomogeneous coordinates in (1.11), which will play an important role in the later $\EUC$-visualization of the prism tilings in $\SLR$,
are given by
\begin{equation}
\begin{gathered}
x=\frac{x^1}{x^0}=\tan{\phi}, \\
y=\frac{x^2}{x^0}=\tanh{r} \frac{\cos{(\theta-\phi)}}{\cos{\phi}}, \\
z=\frac{x^3}{x^0}=\tanh{r} \frac{\sin{(\theta-\phi)}}{\cos{\phi}}.
\end{gathered} \tag{1.11}
\end{equation}

The infinitesimal arc-length-square can be derived by the standard pull back method. 
By $T^{-1}$-action of (1.6) on the differentials $(\mathrm{d}x^0;\mathrm{d}x^1;\mathrm{d}x^2;\mathrm{d}x^3)$, we obtain 
that in this parametrization
the infinitesimal arc-length-square
at any point of $\SLR$ is the following:
\begin{equation}
   \begin{gathered}
      (\mathrm{d}s)^2=(\mathrm{d}r)^2+\cosh^2{r} \sinh^2{r}(\mathrm{d}\theta)^2+\big[(\mathrm{d}\phi)+\sinh^2{r}(\mathrm{d}\theta)\big]^2.
       \end{gathered} \tag{1.12}
     \end{equation}
Hence we get the symmetric metric tensor field $g_{ij}$ on $\SLR$ by components:
     \begin{equation}
       g_{ij}:=
       \begin{pmatrix}
         1&0&0 \\
         0&\sinh^2{r}(\sinh^2{r}+\cosh^2{r})& \sinh^2{r} \\
         0&\sinh^2{r}&1 \\
         \end{pmatrix}, \tag{1.13}
     \end{equation}
and
\begin{equation}
     \mathrm{d}V=\sqrt{\det(g_{ij})}~dr ~\mathrm{d} \theta ~\mathrm{d} \phi= \frac{1}{2}\sinh(2r) \mathrm{d}r ~\mathrm{d} \theta~ \mathrm{d} \phi
     \notag
     \end{equation}
as the volume element in hyperboloid coordinates.
The geodesic curves of $\SLR$ are generally defined as having locally minimal arc length between any two of their (close enough) points.

By (1.13) the second order differential equation system of the $\SLR$ geodesic curve is the following:
\begin{equation}
\begin{gathered}
\ddot{r}=\sinh(2r)~\! \dot{\theta}~\! \dot{\phi}+\frac12 \big( \sinh(4r)-\sinh(2r) \big)\dot{\theta} ~\! \dot{\theta},\\
\ddot{\phi}=2\dot{r}\tanh{(r)}(2\sinh^2{(r)}~\! \dot{\theta}+ \dot{\phi}),\\ \ddot{\theta}=\frac{2\dot{r}}{\sinh{(2r)}}\big((3 \cosh{(2r)}-1)
\dot{\theta}+2\dot{\phi} \big). \tag{1.14}
\end{gathered}
\end{equation}
We can assume, by the homogeneity, that the starting point of a geodesic curve is the origin $(1,0,0,0)$.  
Moreover, $r(0)=0,~ \phi(0)=0,~ \theta(0)=0,~ \dot{r}(0)=\cos(\alpha),~ \dot{\phi}(0)=\sin(\alpha)=-\dot{\theta}(0)$ are the initial values
in Table 1 for the solution of (1.14),
and so the unit velocity will be achieved.
\smallbreak
\centerline{\vbox{
\halign{\strut\vrule~\hfil $#$ \hfil~\vrule
&\quad \hfil $#$ \hfil~\vrule
&\quad \hfil $#$ \hfil\quad\vrule
&\quad \hfil $#$ \hfil\quad\vrule
&\quad \hfil $#$ \hfil\quad\vrule
\cr
\noalign{\vskip2pt}
\noalign{\hrule}
\multispan2{\strut\vrule\hfill{\bf  Table 1} \hfill\vrule}%
\cr
\noalign{\hrule}
\noalign{\vskip2pt}
{\rm Types} & {}  \cr
\noalign{\vskip2pt}
\noalign{\hrule}
\noalign{\vskip2pt}
\begin{gathered} 0 \le \alpha < \frac{\pi}{4} \\ (\bH^2-{\rm like~direction}) \end{gathered}
& \begin{gathered}  r(s,\alpha)={\mathrm{arsinh}} \Big( \frac{\cos{\alpha}}{\sqrt{\cos{2\alpha}}}\sinh(s\sqrt{\cos{2\alpha}}) \Big) \\
\theta(s,\alpha)=-{\mathrm{arctan}} \Big( \frac{\sin{\alpha}}{\sqrt{\cos{2\alpha}}}\tanh(s\sqrt{\cos{2\alpha}}) \Big) \\
\phi(s,\alpha)=2\sin{\alpha} s + \theta(s,\alpha) \end{gathered} \cr
\noalign{\vskip2pt}
\noalign{\hrule}
\noalign{\vskip2pt}
\begin{gathered} \alpha=\frac{\pi}{4} \\ ({\rm light~direction}) \end{gathered} &
\begin{gathered}  r(s,\alpha)={\mathrm{arsinh}} \Big( \frac{\sqrt{2}}{2} s \Big) \\
\theta(s,\alpha)=-{\mathrm{arctan}} \Big( \frac{\sqrt{2}}{2} s \Big) \\
\phi(s,\alpha)=\sqrt{2} s +\theta(s,\alpha) \end{gathered} \cr
\noalign{\vskip1pt}
\noalign{\hrule}
\noalign{\vskip2pt}
\begin{gathered} \frac{\pi}{4}  < \alpha \le \frac{\pi}{2} \\ ({\rm fibre-like~direction}) \end{gathered} &
\begin{gathered}  r(s,\alpha)={\mathrm{arsinh}} \Big( \frac{\cos{\alpha}}{\sqrt{-\cos{2\alpha}}}\sin(s\sqrt{-\cos{2\alpha}}) \Big) \\
\theta(s,\alpha)=-{\mathrm{arctan}} \Big( \frac{\sin{\alpha}}{\sqrt{-\cos{2\alpha}}}\tan(s\sqrt{-\cos{2\alpha}}) \Big) \\
\phi(s,\alpha)=2\sin{\alpha} s + \theta(s,\alpha) \end{gathered}  \cr
\noalign{\vskip2pt}
\noalign{\hrule}
\noalign{\hrule}}}}
\smallbreak
The equation of the geodesic curve in the hyperboloid model has been determined in \cite{DESS}, with the usual geographical sphere coordiantes
$(\lambda, \alpha)$, as longitude and altitude, respectively, from the general starting position of (1.10), (1.11),
$(-\pi < \lambda \le \pi, ~ -\frac{\pi}{2}\le \alpha \le \frac{\pi}{2})$,
and the arc-length parameter $0 \le s \in \bR$. The Euclidean coordinates $X(s,\lambda,\alpha)$,
$Y(s,\lambda,\alpha)$, $Z(s,\lambda,\alpha)$ of the geodesic curves can be determined by substituting the results of Table 1 (see \cite{DESS}) into the
equations (1.10) and (1.11) as follows
\begin{equation}
\begin{gathered}
X(s,\lambda,\alpha)=\tan{(\phi(s,\alpha))}, \\
Y(s,\lambda,\alpha)=\frac{\tanh{(r(s,\alpha))}}{\cos{(\phi(s,\alpha)}} \cos \big[ \theta(s,\alpha)-\phi(s,\alpha)+\lambda \big],\\
Z(s,\lambda,\alpha)=\frac{\tanh{(r(s,\alpha))}}{\cos{(\phi(s,\alpha)}} \sin \big[ \theta(s,\alpha)-\phi(s,\alpha)+\lambda \big].
\end{gathered} \tag{1.15}
\end{equation}
\section{Geodesic balls in $\SLR$}
\begin{Definition}
The {\rm distance} $d(P_1,P_2)$ between the points $P_1$ and $P_2$ is defined by the arc length of the geodesic curve
from $P_1$ to $P_2$.
\end{Definition}
The numerical approximation of the distance $d(O,P)$, by Table 1 and (1.15) for given $P(X,Y,Z)$ from the origin $O$, will not be detailed here.
\begin{Definition}
The {\rm geodesic sphere} of radius $\rho$ (denoted by $S_{P_1}(\rho)$) with the center in point $P_1$ is defined as the set of all points
$P_2$ with the condition $d(P_1,P_2)=\rho$. Moreover, we require that the geodesic sphere is a simply connected
surface without selfintersection.
\end{Definition}
\begin{Definition}
The body of the geodesic sphere of centre $P_1$ and with radius $\rho$ is called {\rm geodesic ball}, denoted by $B_{P_1}(\rho)$,
i.e., $Q \in B_{P_1}(\rho)$ iff $0 \leq d(P_1,Q) \leq \rho$.
\end{Definition}
Fig.~2.a shows a geodesic sphere of radius $\rho=1.3$ with centre $O$ and Fig.~2.b shows its intersection with the $(x,z)$ plane.
\begin{figure}[ht]
\centering
\includegraphics[width=7cm]{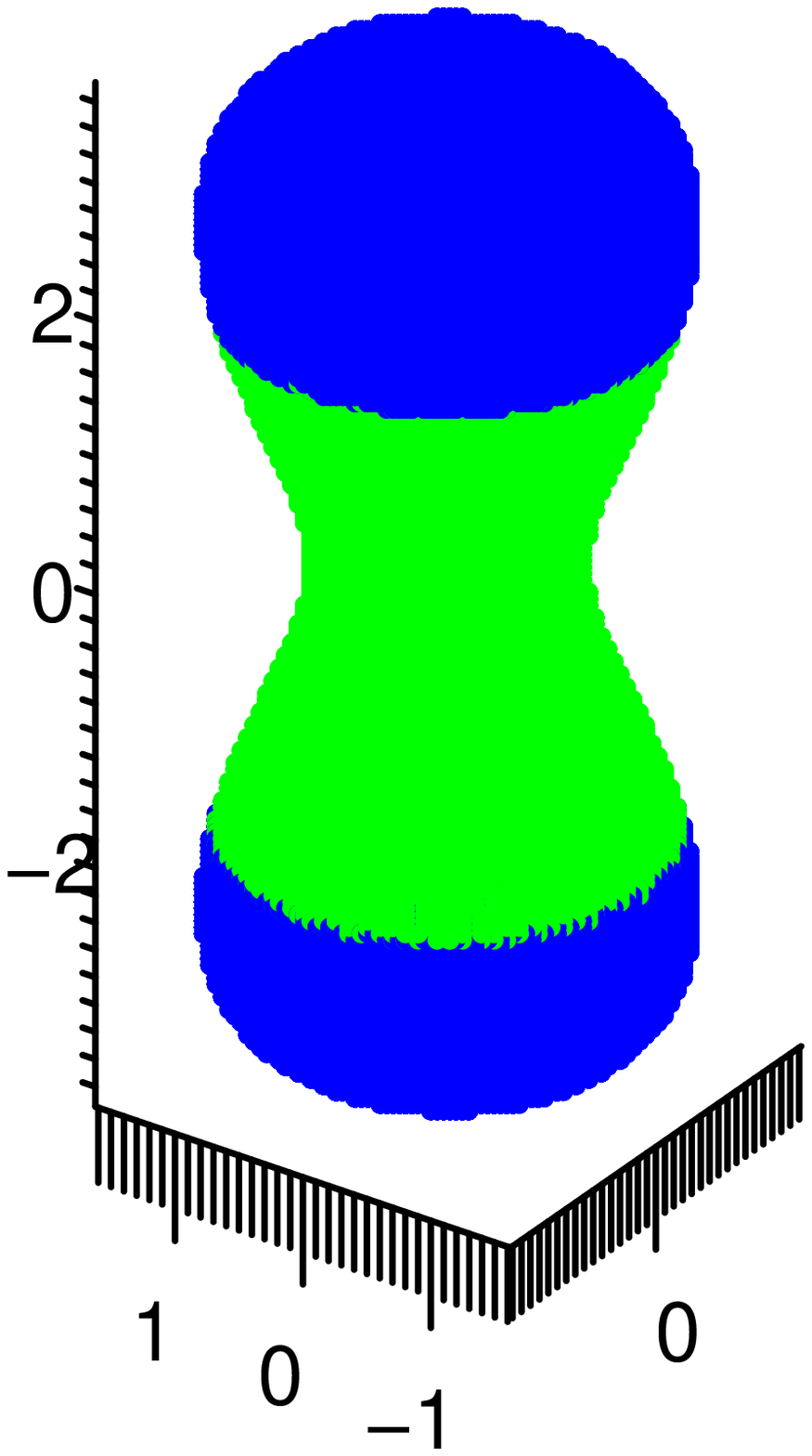} \includegraphics[width=6cm]{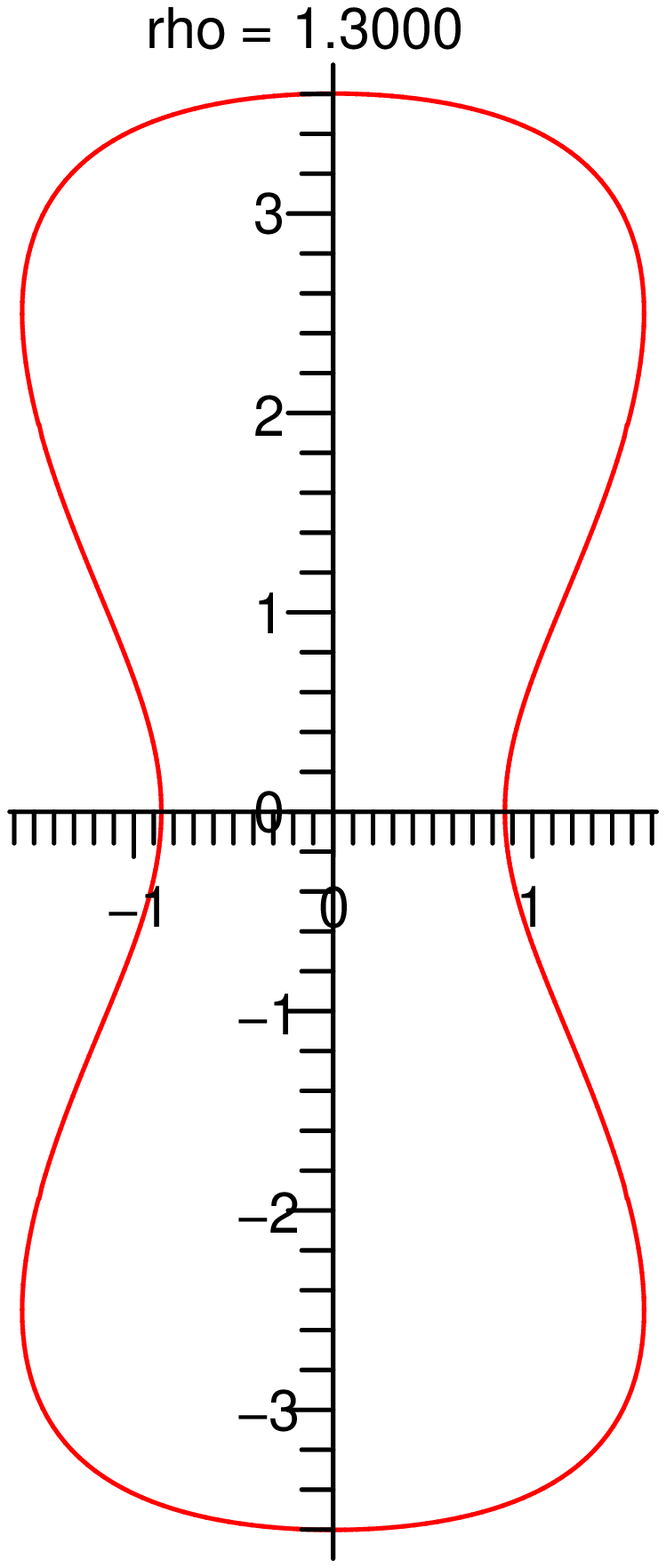}
\caption{a,~b}
\label{}
\end{figure}
 From (1.15) it follows that $S(\rho)$ is a simply connected surface in $\mathbf{E}^3$ and $\SLR$, respectively, if $\rho \in [0,\frac{\pi}{2})$.
 If $\rho\ge \frac{\pi}{2}$ then the universal cover should be discussed.
 {\it Therefore, we consider geodesic spheres and balls only with radii $\rho \in [0,\frac{\pi}{2})$ in the following}. 
 These will be satisfactory for our cases.
\subsection{The volume of a geodesic ball}

The volume formula of the geodesic ball $B(\rho)$ follows from the metric tensor $g_{ij}$.
We obtain the connection between the hyperboloid coordinates $(r,\theta,\phi)$ and the geographical coordinates $(s,\lambda,\alpha)$
by Table 1 and by (1.15). Therefore, the volume of the geodesic ball of radius $\rho$ can be computed by the
following
\begin{theorem}
\begin{equation}
\begin{gathered}
Vol(B(\rho))=\int_{B} \frac{1}{2}\sinh(2r)~ dr ~d \theta~ d \phi  = \\ = 4\pi \int_{0}^{\rho} \int_{0}^{\frac{\pi}{4}}
\frac{1}{2}\sinh(2r(s,  \alpha ))\dot |J_1|  ~ d \alpha \ ds 
\\ +4 \pi \int_{0}^{\rho} \int_{\frac{\pi}{4}}^{\frac{\pi}{2}}
\frac{1}{2}\sinh(2r(s,  \alpha ))\dot |J_2|  ~ d \alpha \ ds
\tag{2.1}
\end{gathered}
\end{equation}
where $|J_1|=
\left|\begin{array}{cccc} \frac{\partial r}{\partial s} & \frac{\partial r}{\partial \alpha} \\
\frac{\partial \phi}{\partial s} & \frac{\partial \phi}{\partial \alpha} \end{array} \right|$
and similarly $|J_2|$ ~ (by Table 1 and  $\frac{\partial \theta}{\partial \lambda}=1$) are the corresponding Jacobians.
\end{theorem}
\begin{figure}[ht]
\centering
\includegraphics[width=8cm]{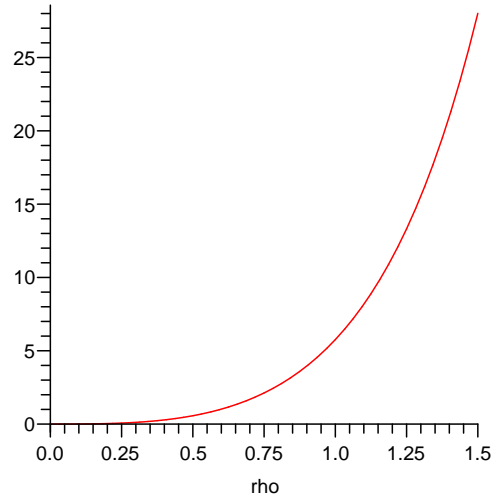}
\caption{The increasing function $\rho \mapsto Vol(B(\rho))$. }
\label{}
\end{figure}
The complicated formulas above need numerical approximations by computer (see Fig.~3).
\section{Regular prism tilings and their space groups $\mathbf{pq2_1}$}
In \cite{Sz13-1} we have defined and described the regular prisms and prism tilings with a space group class $\Gamma=\mathbf{pq2_1}$ of $\SLR$. 
These will be summarized in this section.
\begin{Definition}
Let $\cP^i$ be an infinite solid that is bounded by certain surfaces that can be
determined (in \cite{Sz13-1}) by ,,side fibre lines" passing through the
vertices of a regular $p$-gon $\cP^b$ lying in the base plane.
The images of solids $\cP^i$ by $\SLR$ isometries are called {\rm infinite regular $p$-sided prisms}.
Here regular means that the side surfaces are congruent to each other under rotations about a fiber
line (e.g. through the origin).
\end{Definition}
The common part of $\cP^i$ with the base plane is the {\it base figure} of $\cP^i$ that is denoted by $\cP$ and its vertices coincide
with the vertices of $\cP^b$, {\bf but $\cP$ is not assumed to be a polygon}.
\begin{Definition}
A {\rm bounded regular $p$-sided prism} is analogously defined if the face of the base figure $\cP$ and its translated copy $\cP^t$
by a fibre translation by (1.3) and so (1.5) are also introduced.
The faces $\cP$ and $\cP^t$ are called {\rm cover faces}.
\end{Definition}
\begin{Remark}
All cross-sections of a prism generated by fibre translations from the base plane
are congruent. Prisms are named for their base,
e.g. a prism with a pentagonal base is called {\rm pentagonal prism}.
\end{Remark}
We consider regular prism tilings $\cT_p(q)$ by prisms $\cP_p(q)$ where $q$ pieces regularly meet
at each side edge by $q$-rotation.

The following theorem have been proved in \cite{Sz13-1}:
\begin{theorem}
There exist regular infinite prism tilings $\cT_p^i(q)$ in $\SLR$ for each $3 \le p \in \mathbb{N}$ where $\frac{2p}{p-2} < q \in \mathbb{N}$.
For bounded prisms, these are not face-to-face.
\end{theorem}

We assume that the prism $\cP_p(q)$ is a {\it topological polyhedron} having at each vertex
one $p$-gonal cover face (it is not a polygon at all) and two {\it skew quadrangles} which lie on certain side surfaces in the model.
Let $\cP_p(q)$ be one of the tiles of $\cT_p(q)$, $\cP^b$ is centered in the origin with vertices $A_1A_2A_3 \dots A_p$ in the base plane (Fig.~4).
It is clear that the side curves $c_{A_iA_{i+1}}$ $(i=1\dots p, ~ A_{p+1} \equiv A_1)$
of the base figure are derived from each other by $\frac{2\pi}{p}$ rotation about the vertical $x$ axis, so there are congruent in $\SLR$ sense.
The corresponding vertices $B_1B_2B_3 \dots B_p$ are generated by a fibre translation $\tau$ given by (1.3)
with parameter $0< \Phi\in \mathbb{R}$.
\begin{figure}[ht]
\centering
\includegraphics[width=12cm]{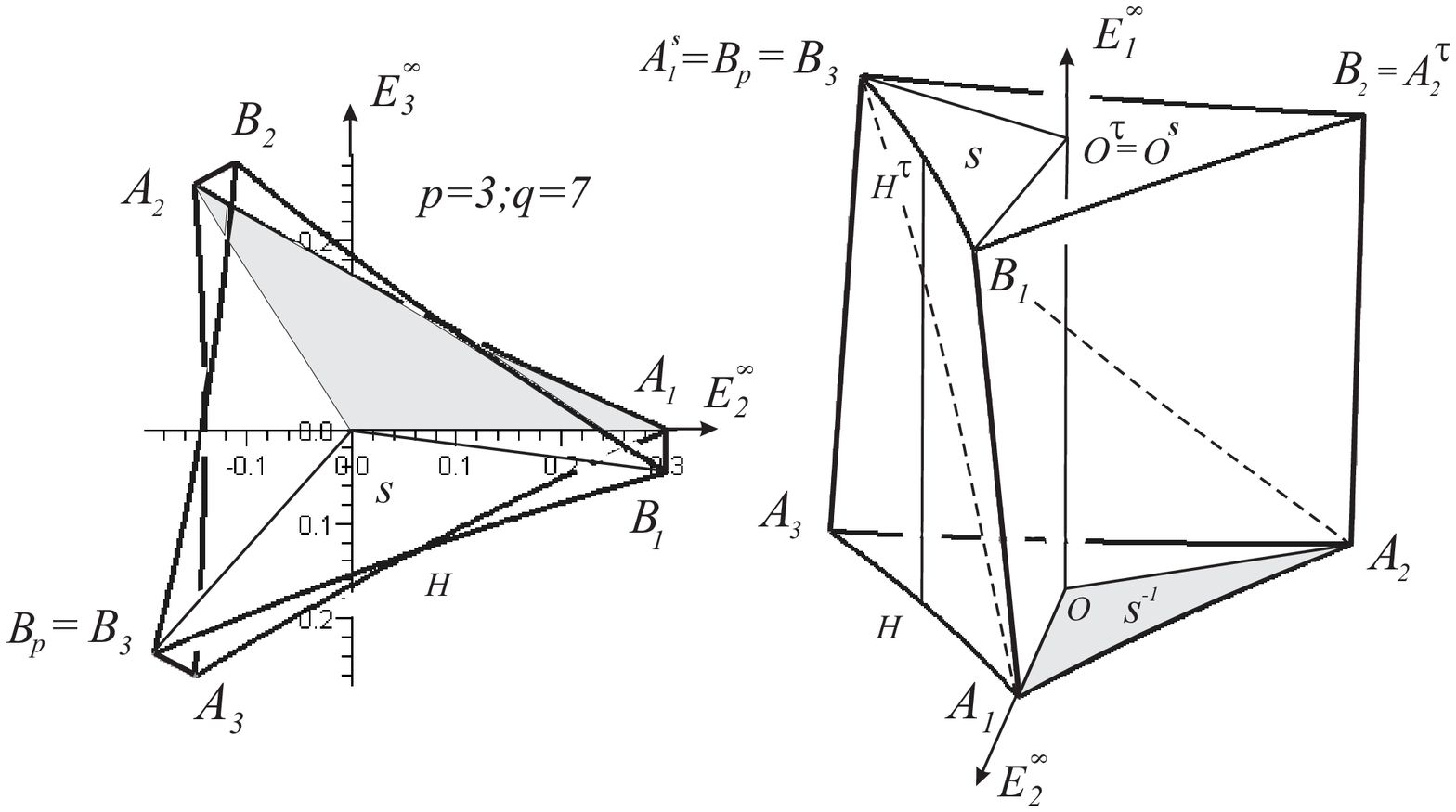}
\caption{}
\label{}
\end{figure}
The fibre lines through the vertices
$A_iB_i$ are denoted by $f_i, \ (i=1, \dots, p)$ and the fibre line through the "midpoint" $H$ of the curve $c_{A_1A_{p}}$ is denoted by
$f_0$. This $f_0$ will be a half-screw axis as follows below.

The tiling $\cT_p(q)$ is generated by a
discrete isometry group $\Gamma_p(q)=\mathbf{pq2_1}$ $\subset Isom(\SLR)$
which is given by its fundamental domain $A_1A_2O A_1^{\bs} A_2^{\bs} O^{\bs}$ a {\it topological polyhedron} and the group presentation
(see Fig.~4 for $p=3$ and \cite{Sz13-1} for details):
\begin{equation}
\begin{gathered}
\mathbf{pq2_1}=\{ \ba,\bb,\bs: \ba^p=\bb^q=\ba \bs \ba^{-1} \bs^{-1}= \bb \ba \bb \bs^{-1}=\mathbf{1} \}= \\
= \{ \ba,\bb: \ba^p=\bb^q=\ba \bb \ba \bb \ba^{-1} \bb^{-1} \ba^{-1} \bb^{-1}=\mathbf{1} \}. \tag{3.1}
\end{gathered}
\end{equation}
Here $\ba$ is a ${p}$-rotation about the fibre line through the origin ($x$ axis),  $\bb$ is a ${q}$-rotation about the fibre line trough
$A_1$ and $\bs=\bb \ba \bb$ is a screw motion~ $\bs:~ OA_1A_2 \rightarrow O^{\bs} B_p B_1$. All these can be obtained by formulas (1.7) and (1.8).
Then we get the second presentation in (3.1), i.e.
$\ba\bb\ba\bb=\bb\ba\bb\ba=:\tau$ is a fibre translation. Then $\ba \bb$ is a $\mathbf{2_1}$ half-screw motion about
$f_0=HH^{\tau}$ (look at Fig.~4) that also determines the fibre tarnslation $\tau$ above. This group in (3.1) surprisingly occurred in \S ~ 6 of our paper \cite{MSzV} at double links
$K_{p,q}$.
The coordinates of the vertices $A_1A_2A_3 \dots A_p$ of the base figure and the corresponding vertices $B_1B_2B_3 \dots B_p$
of the cover face can be computed for all given parameters $p,q$ by
\begin{equation}
\tanh(OA_1)=b:=\sqrt{\frac{1-\tan{\frac{\pi}{p}} \tan{\frac{\pi}{q}}} {1+\tan{\frac{\pi}{q}} \tan{\frac{\pi}{q}}}}. \tag{3.2}
\end{equation}
Moreover, the equation of the curve $c_{A_1A_2}$ can be determined as the foot points (see (1.4) and (1.9)) of the corresponding fibre lines.
For example, the data of $\cP_3(q)$ for some $6<q \in \mathbb{N}$ are collected in Table 2 by Maple computations.
\medbreak
\centerline{\vbox{
\halign{\strut\vrule\quad \hfil $#$ \hfil\quad\vrule
&\quad \hfil $#$ \hfil\quad\vrule
\cr
\noalign{\hrule}
\multispan2{\strut\vrule\hfill\bf Table 2 \hfill\vrule}%
\cr
\noalign{\hrule}
\noalign{\vskip2pt}
\noalign{\hrule}
(p, q)& b  \cr
\noalign{\hrule}
(3,7) & \approx 0.30007426  \cr
\noalign{\hrule}
(3,8) & \approx 0.40561640  \cr
\noalign{\hrule}
(3,9) & \approx 0.47611091  \cr
\noalign{\hrule}
(3,10) & \approx 0.50289355 \cr
\noalign{\hrule}
(3,50) & \approx 0.89636657  \cr
\noalign{\hrule}
(3,1000) & \approx 0.99457331  \cr
\noalign{\hrule}
(3,\infty) & 1  \cr
\noalign{\hrule}
}}}
\medbreak
\subsection{The volume of the bounded regular prism $\cP_p(q)$}
The volume formula of a {\it sector-like} 3-dimensional domain $Vol(D(\Phi))$ can be computed by the metric tensor $g_{ij}$ (1.13)
in hyperboloid coordinates. This defined by the base figure $D$ $(=s^{-1})$ lying in the base plane (see Fig.~4) and by
fibre translation $\tau$ given by (1.3) with the height parameter $\Phi=\pi-\frac{2\pi}{p}-\frac{2\pi}{q}$.
\begin{theorem}
Suppose we are given a sector-like region $D$ (illustrated in Fig.~4), so a continuous function $r = r(\theta)$
where the radius $r$ depends upon the polar angle $\theta$. The volume of domain $D(\Phi))$ is derived by the following integral:
\begin{equation}
\begin{gathered}
Vol(D(\Phi))=\int_{D}  \frac{1}{2}\sinh(2r(\theta)) {\mathrm{d}}r~ {\mathrm{d}} \theta  ~{\mathrm{d}} \phi = \\ = \int_0^{\Phi} \int_{\theta_1}^{\theta_2} \int_{0}^{r(\theta)}
\frac{1}{2}\sinh(2r(\theta))~ {\mathrm{d}}r~ {\mathrm{d}} \theta \ {\mathrm{d}} \phi =\Phi \int_{\theta_1}^{\theta_2}
\frac{1}{4}(\cosh(2r(\theta))-1)~ {\mathrm{d}} \theta.
\tag{3.3}
\end{gathered}
\end{equation}
\end{theorem}
Let $\cT_p(q)$ be the regular prism tiling above and let
$\cP_p(q)$ be one of its tiles. We get the following
\begin{theorem}
The volume of the bounded regular prism $\cP_p(q)$ \Big($3 \le p \in \mathbb{N} $, $\frac{2p}{p-2} < q \in \mathbb{N}$\Big) can be computed by the following formula:
\begin{equation}
Vol(\cP_p(q))=Vol(D(p,q,\Phi)) \cdot p, \tag{3.4}
\end{equation}
where $Vol(D(p,q,\Phi))$ is the volume of the sector-like 3-dimensional domain that is given by the sector region $OA_1A_2 \subset \cP$ (see Fig.~4) 
and by
$\Phi=A_1B_1=$ $\pi-\frac{2\pi}{p}-\frac{2\pi}{q}$, the $\SLR$ height of the prism, depending on $p,q$.
\end{theorem}
\section{The optimal geodesic ball packings under $\mathbf{pq2_1}$}
Sphere packing problems concern arrangements of non-overlapping equal spheres, rather balls, which fill a space.
Space is the usual three-dimensional Euclidean space. However, ball (sphere) packing problems can be generalized to the other
$3$-dimensional Thurston geometries. But sometimes a difficult problem is -- similarly to the hyperbolic space -- the exact definition of the packing density.
In \cite{Sz13-2} we extended the problem of finding the densest geodesic
ball packing for the other $3$-dimensional homogeneous geometries (Thurston geometries). In this paper we study the problem in $\SLR$
and develop a procedure for regular prism tilings and their above group $\mathbf{pq2_1}$ in (3.1).

Let $\mathcal T_p(q)$ be a regular prism tiling and let
$\mathcal P_p(q)$ be one of its tiles which is given by its base figure $\mathcal P$ that is centered at the origin with vertices $A_1A_2A_3 \dots A_p$
in the base plane of the model (see Fig.~5) and the corresponding vertices $B_1B_2$ $B_3 \dots B_p$ and $C_1C_2C_3 \dots C_p$
are generated by fibre translations $\tau:= \mathbf{abab} = \mathbf{baba}$ and its inverse, given by (1.3) (1.8) and (3.1)  
with parameter $\Phi$ also to the above group $\mathbf{pq2_1}$.
It can be assumed by symmetry arguments that the optimal geodesic ball is centered in the origin. 
Denote by $B(E_{0}, \rho)$ the geodesic ball of radius $\rho$ centered in $E_{0}(1;0;0;0)$. 
The volume ${vol}(\mathcal P_p(q))$ is given by the parameters $p$, $q$ 
and $\Phi \ge \rho_{opt}$. The images of $\mathcal P_p(q)$
by the discrete group $\mathbf{pq2_1}$ covers the $\SLR$  space without overlap.
For the density of the packing it is sufficient to relate the volume of the optimal ball
to that of the solid $\mathcal P_p(q)$ (see Definition 3.1). 

We study only one case of the multiply transitive geodesic ball packings where the fundamental domains 
of the $\SLR$ space groups $\mathbf{pq2_1}$ are not prisms. 
Let the fundamental domains be derived by the Dirichlet~--- Voronoi cells (D-V cells) where their centers are images of the origin. 
The volume of the fundamental domain and of the D-V cell is the same, respectively, as in the prism case (for any above $(p, q)$ fixed).

These locally densest geodesic ball packings can be determined for all possible fixed integer parameters $(p,q)$. 
The optimal radius $\rho_{opt}$ is 
\begin{equation}
\rho_{opt} = \min \Big\{ \operatorname{artanh} {(OA_1)}, \ \frac{\Phi}{2}=\frac{\pi}{2}-\frac{\pi}{p}-\frac{\pi}{q}, 
\ \frac{d(O,O^{{\mathbf a} {\mathbf b} })}{2} \Big\} , 
\end{equation}
where $d(O,O^{{\mathbf a} {\mathbf b}})$ is the geodesic distance between $O$ and $O^{{\mathbf a} {\mathbf b}}$ by Definition 2.1. 

The maximal density of the above ball packings can be computed for any possible parameters $p,q$. In Table 3 we have summarized some numerical results.
The best density that we found $\approx 0.567362$ for parameters $p=8, q=10$.

\medbreak
\centerline{\vbox{
\halign{\strut\vrule\quad \hfil $#$ \hfil\quad\vrule
&\quad \hfil $#$ \hfil\quad\vrule &\quad \hfil $#$ \hfil\quad\vrule &\quad \hfil $#$ \hfil\quad\vrule &\quad \hfil $#$ \hfil\quad\vrule
\cr
\noalign{\hrule}
\multispan5{\strut\vrule\hfill\bf Table 3 \hfill\vrule}%
\cr
\noalign{\hrule}
\noalign{\vskip2pt}
\noalign{\hrule}
(p, q)& \rho(K_{opt}) &  {Vol(B_K)} & {Vol(\cP_p(q))} & \delta(K_{opt}) \cr
\noalign{\hrule}
(3,11) & 0.237999 & 0.057543 & 0.169931 & 0.338626 \cr
\noalign{\hrule}
(3,12) & 0.261799 & 0.076892 & 0.205617 & 0.373960 \cr
\noalign{\hrule}
(3,13) & 0.279134 & 0.093489 & 0.238467 & 0.392044 \cr
\noalign{\hrule}
(3,14) & 0.287083 & 0.101857 & 0.268561 & 0.379271 \cr
\noalign{\hrule}
(3,50) & 0.350810 & 0.188371 & 0.636918 & 0.295754 \cr
\noalign{\hrule}
(3,1000) & 0.370822 & 0.223543 & 0.812627& 0.275087 \cr
\noalign{\hrule}
(5,7) & 0.493679 & 0.546132 & 1.218594 & 0.448165 \cr
\noalign{\hrule}
(6,8) & 0.654498 & 1.350812 & 2.570209 & 0.525565 \cr
\noalign{\hrule}
(6,9) & {0.692287} & {1.624770} & {2.924327} & {0.555605} \cr
\noalign{\hrule}
(7,9) & 0.772932 & 2.347696 & 4.181962 & 0.561386 \cr
\noalign{\hrule}
(7,10) & 0.789635 & 2.523909 & 4.568217 & 0.552493 \cr
\noalign{\hrule}
(8,10) & 0.860471 & 3.387783 & 5.971111 & 0.567362 \cr
\noalign{\hrule}
(9,11) & 0.930662 & 4.456867 & 7.887074 & 0.565085 \cr
\noalign{\hrule}
(9,3000) & 1.003711 & 5.838784 & 13.410609 & 0.435385 \cr
\noalign{\hrule}
(20,60) & 1.361357 & 18.712577 & 37.065848 & 0.504847 \cr
\noalign{\hrule}
(20,2000) & 1.387192 & 20.205264 & 39.883121 & 0.506612 \cr
\noalign{\hrule}
}}}
\medbreak
\begin{figure}[ht]
\centering
\includegraphics[width=12cm]{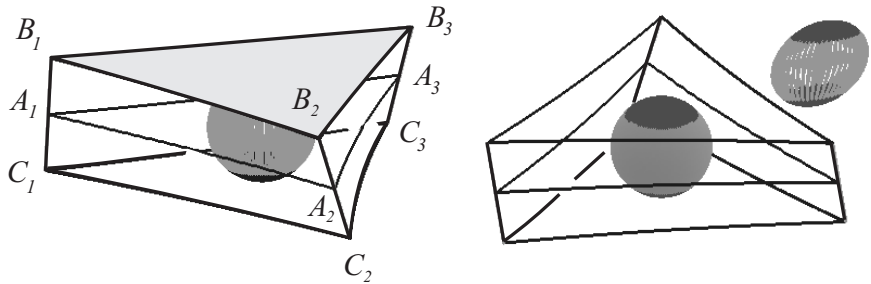} 
\caption{The optimal prism and ball configuration for parameters $p=3$ and $q=7$.}
\label{}
\end{figure}
\begin{Remark}
Surprisingly (at the first glance), the analogous translation ball packings led to larger 
densities, e.g. at $(p, q) = (5, 10000 \rightarrow \infty)$ we obtained the density $0.841700$ close 
to the $\HYP$ upper bound $0.85326$.
\end{Remark}

Our projective method gives us a way of investigating the Thurston geometries. This suits to study and solve similar problems.
(see e.g. \cite{MSz12}, \cite{Sz07}, \cite{Sz12-1}, \cite{Sz13-2}).

%{\bf{Acknowledgement:}}

\end{document}